\documentclass[11pt,oneside]{amsart}
\usepackage{url}
\usepackage[english]{babel}
\usepackage{amsthm, amsmath, amssymb, amsfonts, amscd}
\usepackage[svgnames]{xcolor} 
\usepackage[colorlinks,citecolor=DarkViolet,urlcolor=DarkViolet, bookmarksdepth=3, linkcolor=DarkCyan]{hyperref} 
\usepackage{kpfonts} 
\usepackage{stmaryrd}
\usepackage{listings}
\usepackage{graphicx}
\usepackage{tikz-cd}
\usepackage{indentfirst}
\usepackage{mathrsfs}
\usepackage{accents}
\usepackage{array}
 
\def\sstar{\operatorname{star}}

\def\p{\partial}
\def\k{\Bbbk}

\def\R{\mathbb {R}}
\newtheorem{theorem}{Theorem}

\newtheorem{proposition}[theorem]{Proposition}

\title{Top-heavy Lefschetz elements for regular unimodular embedded lattice simplicial balls with induced boundary}
\author{Mykola Pochekai}
\address{Mykola Pochekai \\ Center for Geometry and Physics \\ Institute for Basic Science \\ Pohang \\ South Korea}
\email{mykola.pochekai@gmail.com}
\begin{document}
\begin{abstract}
We prove that the coordinate Artinian reduction of the Stanley--Reisner ring of a regular unimodular embedded lattice simplicial ball whose boundary is an induced subcomplex admits an explicit top-heavy Lefschetz element.\end{abstract}
\maketitle
\markboth{}{}
\thispagestyle{empty}
\section{Introduction}
There has been recent interest in Lefschetz-type properties for Artinian reductions of Cohen--Macaulay rings and modules arising naturally in combinatorics; see, for instance, \cite{APP21,APP25,PartitionComplex}. Sometimes, one cannot expect the full strong Lefschetz property; the natural replacement is an almost-Lefschetz property and top-heavy Lefschetz elements, see for example \cite[Theorem 46]{PartitionComplex}, \cite[Theorem 6.8]{APP25}, \cite[Theorem 1.1]{Kubitzke}. We prove the following theorem.
\begin{theorem}\label{th:main-theorem}
Let $B$ be a regular unimodular embedded lattice simplicial $d$-ball with induced
boundary, and let $\phi:|B|\to\mathbb R$ be a regularity function. Let
$\mathbb R[B]$ denote the Stanley--Reisner ring of $B$ over $\mathbb R$, and set
$$
A^*(B):=\mathbb R[B]/(\Theta_{\mathrm{coord}}),
$$
where $\Theta_{\mathrm{coord}}$ is the coordinate linear system of parameters
defined in \eqref{eq:lsop-arising-from-parameters}. Then the element
$$
c_\phi=-\sum_{v\in B_0}\phi(v)x_v
$$
is a top-heavy Lefschetz element in $A^*(B)$. Equivalently, multiplication by
$c_\phi^{d-2r}$ induces an injective map
$$
A^r(B)\xrightarrow{\cdot c_\phi^{d-2r}}A^{d-r}(B)
$$
for every $0\le 2r\le d$.
\end{theorem}
The result can be viewed as a variation of the Hard Lefschetz theorem for smooth projective toric varieties \cite[Theorem 12.5.8(e)]{CLS}. Indeed, given a smooth projective fan, the simplicial complex formed by the convex hulls of $0$ and the primitive generators of its cones is a simplicial ball satisfying the hypotheses of Theorem \ref{th:main-theorem}. At the same time, Theorem \ref{th:main-theorem} can be viewed as a coordinate reduction analogue of \cite[Theorem 6.8]{APP25}. The induced boundary assumption implies that $\k[B,\p B]$ is generated in degree $\leq 1$, so \cite[Theorem 6.8]{APP25} establishes the existence of a top-heavy Lefschetz element for generic/suitable Artinian reductions, whereas our theorem establishes this for the coordinate Artinian reduction. 

\subsection*{Acknowledgments.} The author thanks Petter Brändén for insightful discussions and valuable references, and Karim Adiprasito for helpful discussions, references, and suggestions that improved the paper.

\section{Definitions and setup}
By a geometric simplicial complex we mean an abstract finite simplicial complex $B$ together with a map from its vertices to $\R^d$. The map extends affinely to a continuous, $B$-piecewise-linear map $\xi$ from the canonical geometric realization $|B|$ to $\R^d$. If $\sigma$ is a face of $B$, we denote by $\sstar(\sigma)$, $\sstar^\circ(\sigma)$, and $\operatorname{link}(\sigma)$ the closed star, the open star, and the link of $\sigma$ in $B$, respectively. For convenience, we identify simplices, vertices, and subcomplexes of a geometric simplicial complex with their images under $\xi$.
If $B$ is a geometric simplicial complex, we write $B_0$ for the set of all vertices of $B$. We write $\partial B$ for the boundary and $B^\circ$ for the interior of $B$. Abusing notation, we write $B^\circ_0$ for the set of vertices of $B$ that lie strictly in its interior, even though $B^\circ$ is not a simplicial complex. We will also write $\operatorname{link}(\sigma)_0$ and $\sstar(\sigma)_0$ for the set of vertices of the corresponding subcomplexes. For $p_0,\ldots,p_r \in B_0$, we write $p_0 \ldots p_r \in B_r$ if there is an $r$-simplex with vertices $p_0,\ldots,p_r$, and we write $p_0 \ldots p_r \in B^\circ_r$ if, in addition, all vertices $p_0,\ldots,p_r$ lie in the interior of $B$. 
We say that a geometric simplicial complex is nondegenerate if $\xi$ is a topological embedding of $|B|$ into $\R^{d}$ and $B$ is pure of dimension $d$. 
If $B$ is a simplicial complex, we write $\k[B]$ for its Stanley--Reisner ring, graded by total degree, over a field $\k$. For a sequence $\Theta$ of linear forms in $\k[B]$, we write
$$
A^*(B):=\k[B]/(\Theta)
$$
for the corresponding Artinian reduction, the linear system of parameters $\Theta$ will be clear from the context.
We define the linear system of parameters arising from coordinates in the following way:
\begin{equation}\label{eq:lsop-arising-from-parameters}
    \Theta_{\mathrm{coord}} = (\theta_0,\ldots,\theta_d); \qquad \theta_i = \sum_{v \in B_0} \langle v, e_i \rangle x_v \text{ for $1 \leq i \leq d$;} \qquad \theta_0 = \sum_{v \in B_0} x_v,
\end{equation}
where $e_i$ is a canonical basis vector of $\R^d$ and $\langle v, w\rangle$ is a standard scalar product in $\R^d$. We have $\Theta_{\mathrm{coord}} \subset \k[B]$. We also define
\begin{equation}\label{eq:lsop-generic}
    \Theta_{\mathrm{gen}} = (\theta_0,\ldots,\theta_d); \qquad \theta_i = \sum_{v \in B_0} \theta_{iv} x_v \text{ for $0 \leq i \leq d$.} 
\end{equation}
Here $\theta_{iv}$ are formal variables in the field $\tilde \k :=\k(\theta_{iv})_{0 \leq i \leq d, v \in B_0}$, and $\Theta_{\mathrm{gen}} \subset \tilde \k[B]$.

An embedded lattice simplicial $d$-ball $B$ is a nondegenerate geometric simplicial complex such that all vertices have integral coordinates and $|B| \subset \R^d$ is homeomorphic to the standard closed $d$-ball. Such a simplicial ball is unimodular if each maximal simplex is unimodular, that is, has minimal possible volume $1/d!$ among lattice $d$-simplices. We say that an embedded lattice simplicial ball $B$ is regular if there exists a function $\phi:B_0\to\mathbb R$ such that $B$ is obtained as the projection of the lower faces of the lifted point configuration
$$
\{(v,\phi(v)):v\in B_0\}\subset \mathbb R^{d+1}.
$$
The piecewise affine extension of $\phi$ to $|B|$ will be called a regularity function.
We say that a geometric simplicial complex has boundary as an induced subcomplex if every simplex $\sigma \in B$ whose vertices all lie on the boundary is contained entirely in the boundary. Since $B$ is a simplicial $d$-ball, it is Cohen--Macaulay; that is, the face ring $\k[B]$ is Cohen--Macaulay over every field $\k$, which follows from Reisner's theorem \cite[Corollary 5.3.9]{HerzogBruns}.

By a \emph{strong Lefschetz element} we mean an element $\ell \in A^1(B)$ such that multiplication by $\ell^{d-2r}$ induces an isomorphism
$$
A^{r}(B) \xrightarrow{\cdot \ell^{d-2r}} A^{d-r}(B)
$$
for every integer $0 \le 2r \le d$. By a \emph{top-heavy Lefschetz element} $\ell \in A^1(B)$ we mean an element such that multiplication by $\ell^{d-2r}$ induces an injection
$$
A^{r}(B) \xrightarrow{\cdot \ell^{d-2r}} A^{d-r}(B)
$$
for every $0 \leq 2r \le d$.

\section{The proof}
We fix a linear order $<$ on the set of interior vertices $B^\circ_0$. Using this order, we define the unreduced partition complex restricted to interior vertices over field $\k$ as follows:
$$
0 \to \k[B] \to \prod_{p_0 \in B^\circ_0} \k[\sstar(p_0)] \to \ldots \to
\prod_{\substack{p_0 \ldots p_d \in B^\circ_d \\ p_0 < \ldots < p_d}} \k[\sstar(p_0 \ldots p_d)] \to 0,
$$
with the usual alternating restriction maps. The terms are Stanley--Reisner rings of closed stars, while the signs are the Čech signs associated with the corresponding open-star cover. We name this complex $P_{int}^*(B,\k)$, it was introduced in \cite[Subsection 9.2]{PartitionComplex}. The complex is graded in a way that $\k[B]$ has degree $-1$. In \cite[Theorem 50]{PartitionComplex} it is proved that the complex is exact if $B$ is a simplicial ball with induced boundary. We claim that after quotient by a linear system of parameters the complex remains exact.
\begin{proposition} \label{prop:quotient-of-unreduced-partition-complex-is-exact}
    Let $\Delta$ be a simplicial $d$-ball with induced boundary. Let
    $$
    \Theta' = (\theta_0',\ldots,\theta_d') \subset \k[\Delta]_1
    $$
    be a sequence of linear forms such that for all facets $F$ of $\Delta$,
    $$
    \k[F]/(\Theta')\cong \k.
    $$
    Then the quotient complex
    $$
    P^*_{\mathrm{int}}(\Delta,\k)/(\Theta')
    $$
    is exact. In particular $P^*_{\mathrm{int}}/(\Theta_{\mathrm{coord}})$ (in case when $\Delta$ is, additionally, an embedded unimodular simplicial $d$-ball and $\k=\R$) and $P^*_{\mathrm{int}}/(\Theta_{\mathrm{gen}})$ are exact.
\end{proposition}
\begin{proof}
By \cite[Theorem 50]{PartitionComplex}, the complex
$P^*_{\mathrm{int}}$ is exact. We use  \cite[Proposition~1.1.5]{HerzogBruns}: if
$$
N_\bullet:\quad \cdots \to N_m\to N_{m-1}\to\cdots\to N_0\to N_{-1}\to 0
$$
is an exact complex of $R$-modules and $x\in R$ is weakly $N_i$-regular for
every $i$, which means that multiplication by $\cdot x : N_i \to N_i$ is injective, then
$$
N_\bullet\otimes_R R/(x)
$$
is again exact. Applying this
successively, it is enough to show that
$$
\Theta'=(\theta_0',\ldots,\theta_d')
$$
is a regular sequence on every term $P^i_{int}$.
The $i$-th term of $P^*_{\mathrm{int}}$ is a finite product of
rings of the form
$$
\k[\sstar(\tau)],
\qquad
\tau=p_0\cdots p_i\in \Delta_i^\circ.
$$
Thus it is enough to show that $\Theta'$ is a regular sequence on each
$\k[\sstar(\tau)]$. Since $\Delta$ is a simplicial $d$-ball, the closed star
$\sstar(\tau)$ is Cohen--Macaulay of dimension $d$ by Reisner's criterion; see
\cite[Corollary~5.3.9]{HerzogBruns}. Moreover, every facet of $\sstar(\tau)$ is a facet of $\Delta$. Therefore, by the assumption, for every facet $F$ of $\sstar(\tau)$ we have
$$
\k[F]/(\Theta')\cong \k.
$$
By \cite[Theorem~5.1.16(a)]{HerzogBruns}, the sequence $\Theta'$ is a linear system of parameters for $\k[\sstar(\tau)]$. Since $\k[\sstar(\tau)]$ is Cohen--Macaulay, every linear system of parameters
is a regular sequence. Hence $\Theta'$ is regular on each factor
$\k[\sstar(\tau)]$, and therefore on every finite product
$P^i_{int}$. Applying
\cite[Proposition~1.1.5]{HerzogBruns} successively to
$\theta_0',\ldots,\theta_d'$, we conclude that
$$
P^*_{\mathrm{int}}/(\Theta')
$$
is exact.
For $\Theta_{\mathrm{gen}}$, the facet condition holds by genericity. For
$\Theta_{\mathrm{coord}}$, assume that $\Delta$ is an embedded unimodular
simplicial $d$-ball. If $F=\{v_0,\ldots,v_d\}$ is a facet, then the affine
coordinate matrix
$$
\begin{pmatrix}
1 & \cdots & 1\\
v_0 & \cdots & v_d
\end{pmatrix}
$$
has determinant $\pm 1$. Hence the restrictions of
$\theta_0,\ldots,\theta_d$ to $\k[F]$ are linearly independent, and therefore
$\k[F]/(\Theta_{\mathrm{coord}})\cong\k$.
\end{proof}
The point of Proposition~\ref{prop:quotient-of-unreduced-partition-complex-is-exact} is that the partition complex remains exact after quotienting by any sequence of linear forms satisfying the facet condition. For an interior vertex $p\in B_0^\circ$, define a fan $\Sigma_p$ as follows.
We translate $p$ to the origin and, for every face
$\tau\in \operatorname{link}(p)$, set
$$
\sigma_\tau
=
\operatorname{cone}\{v-p : v\in \tau\}
\subset \mathbb R^d.
$$
Then
$$
\Sigma_p
=
\{\sigma_\tau : \tau\in \operatorname{link}(p)\}
$$
is the fan over the link of $p$. Since $p$ is an interior vertex of $B$, the fan $\Sigma_p$ is a complete fan in $\mathbb R^d$. If $B$ is unimodular, then $\Sigma_p$ is smooth. If $B$ is regular, then the
restriction of a regularity function to $\sstar(p)$, after subtracting its
affine part at $p$, gives a strictly convex support function on $\Sigma_p$.
Hence $\Sigma_p$ is projective. We denote by
$$
X_{\sstar(p)}:=X_{\Sigma_p}
$$
the smooth projective toric variety associated with this fan. Now we are ready to prove the Theorem \ref{th:main-theorem}.
\begin{proof}[Proof of Theorem~\ref{th:main-theorem}]
Set
$$
c=c_\phi=-\sum_{v\in B_0}\phi(v)x_v.
$$
By Proposition~\ref{prop:quotient-of-unreduced-partition-complex-is-exact} the quotient $A^*(B):= \R[B]/(\Theta_{\mathrm{coord}})$ by $\Theta_{\mathrm{coord}}:=(\theta_0,\ldots,\theta_d)$, defined by Equation \eqref{eq:lsop-arising-from-parameters}, yields another exact complex
$$0 \to A^*(B) \to \prod_{p_0 \in B_0^\circ} A^*(\sstar(p_0)) \to 
\prod_{\substack{p_0p_1 \in B_1^\circ\\ p_0 < p_1}} A^*(\sstar(p_0p_1)) \to \ldots \to
\prod_{\substack{p_0 \ldots p_d \in B_d^\circ\\ p_0 < \ldots < p_d}} A^*(\sstar(p_0 \ldots p_d)) \to 0,$$
and in particular we obtain an injection
\begin{equation}\label{eq:interior-morphism}
A^*(B) \longrightarrow \prod_{p_0 \in B_0^\circ} A^*(\sstar(p_0)).
\end{equation}
By the Danilov--Jurkiewicz presentation \cite[Theorem~12.4.4]{CLS}, every $A^*(\sstar(p))$ is isomorphic to the cohomology ring
$$H^{2*}(X_{\sstar(p)},\R)$$ 
of the smooth projective toric variety $X_{\sstar(p)}$. By the toric ampleness criterion \cite[Theorem~6.1.14]{CLS}, the restriction 
$c|_{\sstar(p)} := - \sum_{v \in \operatorname{link}(p)_0} (\phi(v)-\phi(p))x_v$  
of the element
$$c = - \sum_{v \in B_0} \phi(v) x_v$$
is an ample element, since $v \mapsto \phi(v) - \phi(p)$ is a strictly convex support function on a fan $\Sigma_p$, so by \cite[Theorem~12.5.8(e)]{CLS}, and the discussion right after the theorem, the element 
$c|_{\sstar(p)}= - \sum_{v \in \operatorname{link}(p)_0} (\phi(v)-\phi(p)) x_v$ 
is a strong Lefschetz element in $A^*(\sstar(p))$ .

Suppose $2k \le d$ and $y \in A^{k}(B)$ is a nonzero element. By the injectivity of \eqref{eq:interior-morphism}, the image $y_q$ of $y$ in $A^{k}(\sstar(q))$ is nonzero for some $q \in B^\circ_0$. Since $c$ restricts to a strong Lefschetz element in $A^{*}(\sstar(q))$, we have $(c|_{\sstar(q)})^{d-2k} y_q \ne 0$ in $A^{d-k}(\sstar(q))$. The commutative square
$$
\begin{tikzcd}
A^{k}(B) \arrow[r] \arrow[d, "{\cdot c^{d-2k}}"'] & A^{k}(\sstar(q)) \arrow[d, "{\cdot c|_{\sstar(q)}^{d-2k}}"] \\
A^{d-k}(B) \arrow[r] & A^{d-k}(\sstar(q))
\end{tikzcd}
$$
then shows that $c^{d-2k} y \ne 0$ in $A^{d-k}(B)$. Therefore, multiplication by $c^{d-2k}$ on $A^{k}(B)$ is injective. Hence $c = c_{\phi}$ is a top-heavy Lefschetz element.
\end{proof}
We want to emphasize that, for the proof to work, it is enough for $\phi$ to be
locally regular: for every interior vertex $p \in B_0^{\circ}$, the function
$v \mapsto \phi(v)-\phi(p)$ must define a strictly convex support function on
$\Sigma_p$.

\bibliographystyle{alpha}
\bibliography{Biblio}{}
\end{document}